\documentclass[11pt,reqno]{amsart}
\usepackage{graphicx}
\usepackage{verbatim}
\usepackage{textcomp}
\usepackage{amssymb}
\usepackage{cite}
\usepackage{amsmath}
\usepackage{latexsym}
\usepackage{amscd}
\usepackage{amsthm}
\usepackage{mathrsfs}
\usepackage{bm}
\usepackage{url}
\usepackage{hyperref}
\usepackage{bookmark}
\allowdisplaybreaks[3]
\newtheorem{thm}{Theorem}[section]
\newtheorem{corr}[thm]{Corollary}
\newtheorem{lem}[thm]{Lemma}

\theoremstyle{definition}

\theoremstyle{remark}
\newtheorem{rem}{Remark}[section]
\numberwithin{equation}{section}
\begin{document}
\title[Gradient estimates for positive weak solution]
{Gradient estimates for positive weak solution to $\Delta_pu+au^{\sigma}=0$ on Riemannian manifolds}

\author{Guangyue Huang }

\author{Qi Guo }

\author{Lujun Guo }

\address{College of Mathematics and Information Science,
Henan Normal University, Xinxiang 453007, P.R. China}

\email{hgy@htu.edu.cn(G. Huang) }

\email{qiguo2022@126.com(Q. Guo) }

\email{lujunguo@htu.edu.cn(L. Guo)}

\thanks{The research of authors is supported by NSFC (Nos. 11971153, 12126319).}

\begin{abstract}
In this paper, we study gradient estimates for positive weak solutions to the following $p$-Laplacian equation
$$\Delta_pu+au^{\sigma}=0$$
on a Riemannian manifold, where $p>1$ and $a,\sigma$ are two nonzero real constants. By virtue of the Morser iteration technique, we derive some gradient estimates, which show that when the Ricci curvature is nonnegative, the above equation does not admit positive weak solutions under some scopes of $p$.
\end{abstract}



\maketitle

\section{Introduction}
Let $M^{n}$ be an $n$-dimensional complete Riemannian manifold with the dimension $n\geq3$. The study on gradient estimates for the positive harmonic functions to the equation
\begin{align}\label{Int-1}
\Delta u=0
\end{align}
was firstly introduced
by Yau \cite{LY1975}, which states that any positive or bounded solution to the
equation \eqref{Int-1} with the Ricci curvature $R_{ij}\geq0$ must be constant. From then on, the study to related nonlinear equations with respect to the Laplacian (can be traced back to Li \cite{Li1991}) have been paid more attention. For the development in this direction, see \cite{Ma2006,Yang2008,Yang2010,MH2017,HM2015,MHL2018,HHL2013,2HHL2013,HL2014,PW2022,PWW2021,WY2021,WY2023,ZS2020} and the references therein.

In \cite{KN2009}, Kotschwar and Ni studied the following $p$-harmonic functions to the equation
\begin{align}\label{Int-2}
\Delta_p u=0,
\end{align}
and derived a gradient estimate under the assumption that the sectional curvature is bounded from below, where the $p$-Laplacian with respect to $u$ is defined by
$$\Delta_pu:={\rm div}(|\nabla u|^{p-2}\nabla u)\ \ \ {\rm for\ every}\ u\in W^{1,p},$$
which is understood in distribution sense. By using the Morser iteration technique, Wang and Zhang \cite{WZ2011} derived a gradient estimate, where their estimates depend on the Ricci curvature instead of the sectional curvature. Recently, Wang and Wei \cite{WW2023} also use the Morser iteration technique and study the nonexistence of positive weak solution to
\begin{align}\label{Int-33}
\Delta u+au^{q+1}=0
\end{align}
on a Riemannian manifold with $a>0$ and proved that if the Ricci curvature is nonnegative, then the above equation does not admit positive weak solutions provided $q\in (-\infty,\frac{2}{n-1}+\frac{2}{\sqrt{n(n-1)}})$.

In this paper, we continue to study gradient estimates for positive weak solutions to the following $p$-Laplacian equation
\begin{align}\label{Int-3}
\Delta_p u+au^{\sigma}=0,
\end{align}
where $p>1$ and $a,\sigma$ are two nonzero real constants. By virtue of the Morser iteration technique, we derive the following gradient estimate:

\begin{thm}\label{Th-1}
Let $M^{n}$ be an $n$-dimensional complete Riemannian manifold with $R_{ij}\geq-(n-1)K$, where $K$ is a positive constant. Assume that the constants $a,\sigma$ satisfy that either $a>0$ and $\sigma<(p-1)\big[\frac{n+1}{n-1}+\frac{2}{n-1}\sqrt{1-\frac{(p-1)^2}{(n-1)\alpha}}\big]$, or $a<0$ and $\sigma>(p-1)\big[\frac{n+1}{n-1}-\frac{2}{n-1}\sqrt{1-\frac{(p-1)^2}{(n-1)\alpha}}\big]$, where
\begin{equation}\label{alpha}
\alpha=\left\{\begin{array}{l}
\frac{n}{n-1}(p-1)^2, \quad p\in\big(1,3-\frac{2}{n}\big],\\
2(p-1), \quad\quad\ \  p\in(3-\frac{2}{n},2n-1).
\end{array}\right.
\end{equation}
Then for any positive weak solution $u$ to \eqref{Int-3} with $1<p<2n-1$ on the geodesic ball $B_{x_0}(R)$, we have
\begin{align}\label{Int-4}
\frac{|\nabla u|}{u}\leq C_{p,n,\sigma}\Big(\frac{1+\sqrt{K}R}{R}\Big)^\frac{p}{2}\ \ \  on\ B_{x_0}(\frac{R}{2}),
\end{align}
where the constant $C_{p,n,\sigma}$ depends only on $p,\sigma$ and the dimension $n$.
\end{thm}

From the above theorem, the following nonexistence results follow  by taking $R\rightarrow \infty$ immediately:

\begin{corr}\label{corr-4}
Let $M^{n}$ be an $n$-dimensional complete Riemannian manifold with nonnegative Ricci curvature with $1<p<2n-1$ and $\alpha$ given by \eqref{alpha}. Assume that the constants $a,\sigma$ satisfy that either $a>0$ and $\sigma\in\big(-\infty,(p-1)\big[\frac{n+1}{n-1}+\frac{2}{n-1}\sqrt{1-\frac{(p-1)^2}{(n-1)\alpha}}\big]\big)$, or $a<0$ and $\sigma\in\big((p-1)\big[\frac{n+1}{n-1}-\frac{2}{n-1}\sqrt{ 1-\frac{(p-1)^2}{(n-1)\alpha}}\big],+\infty\big)$, then there does not exist any positive weak solution to the equation \eqref{Int-3}.

\end{corr}

On the other hand, by using other methods, we can achieve the following gradient estimate, which does not depends on the upper bound of $p$:

\begin{thm}\label{Th-2}
Let $M^{n}$ be an $n$-dimensional complete Riemannian manifold with $R_{ij}\geq-(n-1)K$, where $K$ is a positive constant. Assume that the constants $a,\sigma$ satisfy that either $a>0$ and $\sigma\leq\frac{n+2}{n}(p-1)$, or $a<0$ and $\sigma\geq\frac{n+2}{n}(p-1)$. Then for any positive weak solution $u$ to \eqref{Int-3} with $p>1$ on the geodesic ball $B_{x_0}(R)$, we have
\begin{align}\label{Int-5}
\frac{|\nabla u|}{u}\leq C_{p,n}\Big(\frac{1+\sqrt{K}R}{R}\Big)^{\frac{p}{2}}\ \ \  on\ B_{x_0}(\frac{R}{2}),
\end{align}
where the constant $C_{p,n}$ depends only on $p$ and the dimension $n$.
\end{thm}	

Similarly, the following nonexistence results also follow  by taking $R\rightarrow \infty$ immediately:

\begin{corr}\label{corr-3}
Let $M^{n}$ be an $n$-dimensional complete Riemannian manifold with nonnegative Ricci curvature. Assume that the constants $a,\sigma$ satisfy that either $a>0$ and $\sigma\leq\frac{n+2}{n}(p-1)$, or $a<0$ and $\sigma\geq\frac{n+2}{n}(p-1)$.
Then there does not exist any positive weak solution to the equation \eqref{Int-3}.
\end{corr}

From the above theorem, the following Harnack inequality follows immediately:

\begin{corr}\label{corr-1}
Let $M^{n}$ be an $n(n>2)$-dimensional complete Riemannian manifold with $R_{ij}\geq-(n-1)K$, where $K$ is a positive constant.
Under the same assumptions as in Theorem \ref{Th-2}, then for any positive weak solution $u$ to \eqref{Int-3}, we have
\begin{align}\label{Int-6}
u(x)\leq e^{C_{p,n,K}}u(y),
\end{align}
where $x,y\in B_{x_0}(\frac{R}{2})$.
\end{corr}

\begin{rem}
In particular, when $R_{ij}\geq0$, the above Harnack inequality \eqref{Int-6} becomes
\begin{align}\label{Int-7}
u(x)\leq e^{C_{p,n}}u(y),
\end{align}
where the uniform constant $C_{p,n}$ does not depend on $R$.
\end{rem}

\begin{rem}
Clearly, the Theorem 1.2 in \cite{WZ2011} of Wang and Zhang follows from by taking $a\rightarrow 0$ in \eqref{Int-5}.
\end{rem}

\begin{rem}
For smooth metric measure spaces, our Theorem \ref{Th-2} still holds under the same conditions. On the other hand, Zhao and Yang \cite{ZY2018} derived similar results under the condition $a>0$ and $\sigma\leq p-1$, while our result holds for
$a>0$ and $\sigma\leq\frac{n+2}{n}(p-1)$, or $a<0$ and $\sigma\geq\frac{n+2}{n}(p-1)$. Obviously, we generalize the results in
\cite{ZY2018} of Zhao and Yang.
\end{rem}

\begin{rem}
When $p=2$ and $a>0$, Theorem 1.4 of Wang and Wei in \cite{WW2023} follows from Corollary \ref{corr-4}. In the case of $a<0$, our results are new. On the other hand, if $p<2n-1$, we can check that
$$\frac{n+2}{n}(p-1)<(p-1)\Big[\frac{n+1}{n-1}+\frac{2}{n-1}\sqrt{1-\frac{(p-1)^2}{(n-1)\alpha}}\Big],$$
which shows that, for $a>0$ and $p\geq2n-1$, our Theorem \ref{Th-2} is new.
\end{rem}

\section{Proof of Theorem \ref{Th-1}}
Let $u$ be a positive weak solution to
\begin{align}\label{Pf-1}
{\Delta}_{p}u+au^{\sigma}=0.
\end{align}
We define
$$v:=(p-1)\log u,$$
which is equivalent to $u=e^\frac{v}{p-1}$. Then, we obtain from \eqref{Pf-1} that $v$ satisfies
\begin{equation}\label{Pf-2}
\Delta_pv=-|\nabla v|^p-a(p-1)^{p-1}e^{\big(\frac{\sigma}{p-1}-1\big)v}.
\end{equation}
For
$$f:=|\nabla v|^p,$$
using the elliptic operator $\mathcal{L}$ defined by
$$\mathcal{L}={\rm div}(|\nabla v|^{p-2}A\nabla\cdot),$$
where $A=Id+(p-2)\frac{\nabla v\otimes\nabla v}{|\nabla v|^2}$, we have the following:

\begin{lem}\label{Le-1}
For any $p\in (1,2n-1)$, at the point where $\{f\neq0\}$, we have
\begin{align}\label{Pf-3}
\mathcal{L}(f)\geqslant&-p(n-1)Kf^\frac{2p-2}{p}+\Big(1-\frac{(p-1)^2}{(n-1)\alpha}\Big)\frac{pa^2h^2}{n-1}
+\frac{p}{n-1}f^2\notag\\&
+\Big(\frac{2(p-1)}{n-1}-p\Big)f^\frac{p-2}{p}\nabla f\nabla v +aph\Big[\frac{2}{n-1}-\Big(\frac{\sigma}{p-1}-1\Big)\Big]f,
\end{align}
where $h=(p-1)^{p-1}e^{\big(\frac{\sigma}{p-1}-1\big)v}$.
\end{lem}

\proof
For $p>1$, the following Bochner formula(see the formula (2.3) in \cite{Wang-Li-2016}) with respect to the elliptic operator $\mathcal{L}$ is well-known:
\begin{align}\label{Pf-4}
\frac{1}{p}\mathcal{L}(|\nabla v|^p)=&|\nabla v|^{2p-4}|{\rm Hess}\,v|^2_A+|\nabla v|^{2p-4}{\rm Ric}(\nabla v,\nabla v)\notag\\
&+|\nabla v|^{p-2}\nabla v\nabla \Delta_{p}v,
\end{align}
where $|{\rm Hess}\,v|^2_A=A^{ik}A^{jl}v_{ij}v_{kl}$. Let $\{e_1,e_2,\cdots,e_n\}$ be an orthonormal frame at the point $x_0\in M$ with $e_1\parallel \nabla v$. Then, at the point $x$, we have $v_1=|\nabla v|$ and $v_2=\cdots v_n=0$. Then, \eqref{Pf-2} can be written as
$$|\nabla v|^{p-2}\Delta v+\frac{p-2}{p}|\nabla v|^{-2}\nabla v\nabla|\nabla v|^p=-|\nabla v|^p-a(p-1)^{p-1}e^{\big(\frac{\sigma}{p-1}-1\big)v},$$
which is equivalent to
\begin{equation}\label{Pf-5}
(p-2)|\nabla v|^{p-2}v_{11}+|\nabla v|^{p-2}\sum^{n}_{i=1}v_{ii}=-|\nabla v|^p-ah.
\end{equation}
Then, it follows from \eqref{Pf-5} that
\begin{align}\label{Pf-6}
(p-1)v_{11}+\sum^{n}_{i=2}v_{ii}=-f^\frac{2}{p}-ahf^\frac{2-p}{p},
\end{align}
which implies
\begin{align}\label{Pf-7}
\sum^{n}_{i,j=1}v_{ij}^2\geq&v^2_{11}+2\sum^{n}_{i=2}v^2_{i1}+\sum^{n}_{i=2}v^2_{ii}\notag\\
\geq&v^2_{11}+2\sum^{n}_{i=2}v^2_{i1}+\frac{1}{n-1}\big(\sum^{n}_{i=2}v_{ii}\big)^2\notag\\
=&v^2_{11}+2\sum^{n}_{i=2}v^2_{i1}+\frac{1}{n-1}\Big[f^\frac{2}{p}+ahf^\frac{2-p}{p}+(p-1)v_{11}\Big]^2\notag\\
=&\Big[1+\frac{(p-1)^2}{n-1}\Big]v^2_{11}+2\sum^{n}_{i=2}v^2_{i1}+\frac{1}{n-1}\Big[f^\frac{4}{p}
+a^2h^2f^\frac{4-2p}{p}+2ahf^\frac{4-p}{p}\notag\\
&+2(p-1)v_{11}f^\frac{2}{p}+2a(p-1)hv_{11}f^\frac{2-p}{p}\Big]\notag\\
=&\Big[1+\frac{(p-1)^2}{n-1}\Big]v^2_{11}+2\sum^{n}_{i=2}v^2_{i1}+\frac{1}{n-1}f^\frac{4}{p}
+\frac{a^2h^2}{n-1}f^\frac{4-2p}{p}\notag\\&+\frac{2ah}{n-1}f^\frac{4-p}{p}+\frac{2(p-1)}{n-1}v_{11}f^\frac{2}{p}
+\frac{2a(p-1)h}{n-1}v_{11}f^\frac{2-p}{p}.
\end{align}
By virtue of $pfv_{11}=\nabla v\nabla f$ and $\sum^{n}_{i=1}v^2_{i1}=\frac{1}{p^2}f^\frac{2-2p}{p}|\nabla f|^2$, we obtain
\begin{align}\label{Pf-8}
|{\rm Hess}\,v|^2_A=&\Big[\delta_{ik}+(p-2)|\nabla v|^{-2}v_iv_k\Big]\Big[\delta_{jl}+(p-2)|\nabla v|^{-2}v_jv_l\Big]v_{ij}v_{kl}\notag\\
=&\sum^{n}_{i,j=1}v^2_{ij}+2(p-2)|\nabla v|^{-2}v_{ij}v_{il}v_jv_l+(p-2)^2|\nabla v|^{-4}v_{ij}v_{kl}v_iv_jv_kv_l\notag\\
=&\sum^{n}_{i,j=1}v^2_{ij}+2(p-2)\sum^{n}_{i=1}v^2_{i1}+(p-2)^2v^2_{11}\notag\\
=&\Big[1+(p-2)^2+\frac{(p-1)^2}{n-1}\Big]v^2_{11}+2(p-2)\sum^{n}_{i=1}v^2_{i1}+2\sum^{n}_{i=2}v^2_{i1}
+\frac{1}{n-1}f^{\frac{4}{p}}\notag\\
&+\frac{a^2h^2}{n-1}f^\frac{4-2p}{p}+\frac{2ah}{n-1}f^\frac{4-p}{p}+\frac{2(p-1)}{n-1}v_{11}f^\frac{2}{p}
+\frac{2a(p-1)h}{n-1}v_{11}f^\frac{2-p}{p}\notag\\
=&2(p-2)\sum^{n}_{i=1}v^2_{i1}+\min\Big\{2,1+(p-2)^2+\frac{(p-1)^2}{n-1}\Big\}   \sum^{n}_{i=1}v^2_{i1}+\frac{1}{n-1}f^\frac{4}{p}\notag\\
&+\frac{a^2h^2}{n-1}f^\frac{4-2p}{p}+\frac{2ah}{n-1}f^\frac{4-p}{p}+\frac{2(p-1)}{n-1}v_{11}f^\frac{2}{p}+\frac{2a(p-1)h}{n-1}v_{11}f^\frac{2-p}{p}\notag\\
\geq&\frac{\alpha}{p^2}|\nabla f|^2f^\frac{2-2p}{p}+\frac{1}{n-1}f^\frac{4}{p}+\frac{a^2h^2}{n-1}f^\frac{4-2p}{p}+\frac{2ah}{n-1}f^\frac{4-p}{p}\notag\\
&+\frac{2(p-1)}{p(n-1)}f^\frac{2-p}{p}\nabla f\nabla v+\frac{2a(p-1)h}{p(n-1)}f^\frac{2-2p}{p}\nabla f\nabla v,
\end{align}
where $\alpha$ is defined by \eqref{alpha}.
According to $\eqref{Pf-4},$ we obtain
\begin{align}\label{Pf-9}
\mathcal{L}(f)=&p|\nabla v|^{2p-4}|{\rm Hess}\,v|^2_A+p|\nabla v|^{2p-4}{\rm Ric}(\nabla v,\nabla v)\notag\\
&+p|\nabla v|^{p-2}\nabla v\nabla\Delta_{p}v\notag\\
\geq&\frac{\alpha}{p}|\nabla f|^2f^{-\frac{2}{p}}+\frac{p}{n-1}f^2+\frac{pa^2h^2}{n-1}+\frac{2aph}{n-1}f\notag\\
&+\Big[\frac{2(p-1)}{n-1}-p\Big] f^\frac{p-2}{p}\nabla f\nabla v-\frac{2|a|(p-1)h}{n-1}f^{-\frac{1}{p}}|\nabla f|\notag\\
&-p(n-1)Kf^\frac{2p-2}{p}-apf^\frac{p-2}{p}\nabla v\nabla\Big[(p-1)^{p-1}e^{\big(\frac{\sigma}{p-1}-1\big)v}\Big]\notag\\
=&\frac{\alpha}{p}|\nabla f|^2f^{-\frac{2}{p}}+\frac{p}{n-1}f^2+\frac{pa^2h^2}{n-1}+aph\Big[\frac{2}{n-1}
-(\frac{\sigma}{p-1}-1)\Big]f\notag\\
&+\Big[\frac{2(p-1)}{n-1}-p\Big]f^\frac{p-2}{p}\nabla f\nabla v-\frac{2|a|(p-1)h}{n-1}f^{-\frac{1}{p}}|\nabla f|\notag\\
&-p(n-1)Kf^\frac{2p-2}{p}.
\end{align}
For the above $\alpha$, when $1<p<2n-1$, we have
$$\frac{(p-1)^2}{(n-1)\alpha}<1.$$
Applying the Cauchy inequality
$$\frac{\alpha}{p}|\nabla f|^2f^{-\frac{2}{p}}+\frac{(p-1)^2}{(n-1)\alpha}\frac{pa^2h^2}{n-1}\geq\frac{2|a|(p-1)h}{n-1}f^{-\frac{1}{p}}|\nabla f|$$
into \eqref{Pf-9} gives the desired estimate \eqref{Pf-3} and the proof of Lemma \ref{Le-1} is completed.
\endproof

Without loss of generality, we may assume that $f>0$ (because of at those points where $f=0$, the estimate \eqref{Int-4} is trivial). For simplicity, we let $B_R:=B_{x_0}(R)$. For any non-negative function $\psi\in W_0^{1,2}(B_R)$, we have
\begin{align}\label{Pf-10}
\int_{B_R}&[f^{1-\frac{2}{p}}\nabla{f}+(p-2)
f^{{1-\frac{4}{p}}}(\nabla{f}\nabla{v})\nabla{v}]\nabla{\psi}+{\frac{p}{n-1}}\int_{B_R} f^2\psi\notag\\\leq&(n-1)Kp\int_\Omega f^{2-\frac{2}{p}}\psi-\Big[\frac{2(p-1)}{n-1}-p\Big]\int_{B_R} f^{1-\frac{2}{p}}\psi\nabla{f}\nabla{v}\notag\\
&+\int_{B_R} \Big\{ \Big[\Big(\frac{\sigma}{p-1}-1\Big)-\frac{2}{n-1}\Big]pah f-\Big[1-\frac{(p-1)^2}{(n-1)\alpha}\Big]\frac{pa^2h^2}{n-1}\Big\}\psi.
\end{align}
In what follows, for the convenience, we let $R_3$ to denote the third term on the right hand sides($L_1$ denotes the first term on the left hand, etc.).

In the case of $a>0$, if $\sigma\leq\frac{n+1}{n-1}(p-1)$, we have
$$\Big[\Big(\frac{\sigma}{p-1}-1\Big)-\frac{2}{n-1}\Big]pah f-\Big[1-\frac{(p-1)^2}{(n-1)\alpha}\Big]\frac{pa^2h^2}{n-1}\leq0.$$
Otherwise, if $\sigma>\frac{n+1}{n-1}(p-1)$, we let
$$B_1=\Big\{ x\in B_R; f\geqslant\frac{\frac{1}{n-1}\big[1-\frac{(p-1)^2}{(n-1)\alpha}\big]}{\frac{\sigma}{p-1}-1-\frac{2}{n-1}}\Big\},$$
$$B_2=\Big\{ x\in B_R; f<\frac{\frac{1}{n-1}\big[1-\frac{(p-1)^2}{(n-1)\alpha}\big]}{\frac{\sigma}{p-1}-1-\frac{2}{n-1}}\Big\}.$$
Then,
\begin{align}\label{Pf-11}
&\int_{B_R} \Bigg\{ \Big[(\frac{\sigma}{p-1}-1)-\frac{2}{n-1}\Big]pah f-\Big[1-\frac{(p-1)^2}{(n-1)\alpha}\Big]\frac{pa^2h^2}{n-1}
\Bigg\}\psi\notag\\
&\leq\int_{B_1}\Bigg\{ \Big[(\frac{\sigma}{p-1}-1)-\frac{2}{n-1}\Big]pah f-\Big[1-\frac{(p-1)^2}{(n-1)\alpha}\Big]\frac{pa^2h^2}{n-1}
\Bigg\}\psi \notag\\
&\leq\frac{p\big[(\frac{\sigma}{p-1}-1)-\frac{2}{n-1}\big]^2}{\frac{4}{n-1}\big[1-\frac{(p-1)^2}{(n-1)\alpha}\big]}
\int_{B_1}f^2\psi\notag\\
&\leq\frac{p\big[(\frac{\sigma}{p-1}-1)-\frac{2}{n-1}\big]^2}{\frac{4}{n-1}\big[1-\frac{(p-1)^2}{(n-1)\alpha}\big]}\int_{B_R}f^2\psi.
\end{align}
Denote by $\sigma_1:=(p-1)\Big[\frac{n+1}{n-1}+\frac{2}{n-1}\sqrt{1-\frac{(p-1)^2}{(n-1)\alpha}}\Big]$.
A direct computation shows that if $\sigma\in\big(\frac{n+1}{n-1}(p-1),\sigma_1\big)$, we have $$\frac{p\big[(\frac{\sigma}{p-1}-1)-\frac{2}{n-1}\big]^2}{\frac{4}{n-1}\big[1-\frac{(p-1)^2}{(n-1)\alpha}\big]}<\frac{p}{n-1}.$$
Therefore, for
\begin{equation}\label{Pf-12}
\beta_1(n,p,\sigma)=\left\{\begin{array}{l}
\frac{p}{n-1}, \quad \quad\ \quad\ \quad \quad\ \quad\ \quad \  \sigma\leq\frac{n+1}{n-1}(p-1),\\
\frac{p}{n-1}-\frac{p\big[(\frac{\sigma}{p-1}-1)-\frac{2}{n-1}\big]^2}{\frac{4}{n-1}\big[1-\frac{(p-1)^2}{(n-1)\alpha}\big]}, \quad \frac{n+1}{n-1}(p-1)<\sigma<\sigma_1,
\end{array}\right.
\end{equation}
we obtain from \eqref{Pf-10} that
\begin{align}\label{Pf-13}
\int_{B_R}&[f^{1-\frac{2}{p}}\nabla{f}+(p-2)
f^{{1-\frac{4}{p}}}(\nabla{f}\nabla{v})\nabla{v}]\nabla{\psi}+\beta_1(n,p,\sigma)\int_{B_R} f^2\psi\notag\\\leq&(n-1)Kp\int_{B_R} f^{2-\frac{2}{p}}\psi-\Big[\frac{2(p-1)}{n-1}-p\Big]\int_{B_R} f^{1-\frac{2}{p}}\psi\nabla{f}\nabla{v}.
\end{align}

In the case of $a<0$, if $\sigma>\frac{n+1}{n-1}(p-1)$, we have
$$\int_{B_R} \Bigg\{\Big[(\frac{\sigma}{p-1}-1)-\frac{2}{n-1}\Big]pah f-\Big[1-\frac{(p-1)^2}{(n-1)\alpha}\Big]\frac{pa^2h^2}{n-1}\Bigg\}\psi\leq0.$$
Otherwise, for $\sigma\leq\frac{n+1}{n-1}(p-1)$, we denote
$$B_3=\Bigg\{ x\in B_R; f\leq\frac{\frac{1}{n-1}\Big[1-\frac{(p-1)^2}{(n-1)\alpha}\Big]}{\frac{\sigma}{p-1}-1-\frac{2}{n-1}}\Bigg\},$$
$$B_4=\Bigg\{ x\in B_R; f>\frac{\frac{1}{n-1}\Big[1-\frac{(p-1)^2}{(n-1)\alpha}\Big]}{\frac{\sigma}{p-1}-1-\frac{2}{n-1}}\Bigg\}.$$
Then,
\begin{align}\label{Pf-14}
&\int_{B_R}\Bigg\{ \Big[\Big(\frac{\sigma}{p-1}-1\Big)-\frac{2}{n-1}\Big]pah f-\Big[1-\frac{(p-1)^2}{(n-1)\alpha}\Big]\frac{pa^2h^2}{n-1}\Bigg\}\psi\notag\\
&\leq\int_{B_3}\Bigg\{ \Big[\Big(\frac{\sigma}{p-1}-1\Big)-\frac{2}{n-1}\Big]pah f-\Big[1-\frac{(p-1)^2}{(n-1)\alpha}\Big]\frac{pa^2h^2}{n-1}\Bigg\}\ \psi\notag\\
&\leq\frac{p\big[\big(\frac{\sigma}{p-1}-1\big)-\frac{2}{n-1}\big]^2}{\frac{4}{n-1}\big[1-\frac{(p-1)^2}{(n-1)\alpha}\big]}
\int_{B_3}f^2\psi\notag\\
&\leq\frac{p\Big[(\frac{\sigma}{p-1}-1)-\frac{2}{n-1}\Big]^2}{\frac{4}{n-1}\Big[1-\frac{(p-1)^2}{(n-1)\alpha}\Big]}
\int_{B_R}f^2\psi.
\end{align}
Denote by $\sigma_2:=(p-1)\Big[\frac{n+1}{n-1}-\frac{2}{n-1}\sqrt{1-\frac{(p-1)^2}{(n-1)\alpha}}\Big]$.
By a direct computation, we have that if $\sigma\in\big(\sigma_2,\frac{n+1}{n-1}(p-1)\big)$, then $$\frac{p\big[\big(\frac{\sigma}{p-1}-1\big)-\frac{2}{n-1}\big]^2}{\frac{4}{n-1}\big[1-\frac{(p-1)^2}{(n-1)\alpha}\big]}<\frac{p}{n-1}.$$
Therefore, for
\begin{equation}\label{Pf-15}
\beta_2(n,p,\sigma)=\left\{\begin{array}{l}
\frac{p}{n-1}, \quad \quad\ \quad\ \quad \quad\ \quad\ \quad \  \sigma>\frac{n+1}{n-1}(p-1),\\
\frac{p}{n-1}-\frac{p\big[(\frac{\sigma}{p-1}-1)-\frac{2}{n-1}\big]^2}{\frac{4}{n-1}\big[1-\frac{(p-1)^2}{(n-1)\alpha}\big]}, \quad \sigma_2<\sigma\leq\frac{n+1}{n-1}(p-1),
\end{array}\right.
\end{equation}
we obtain from \eqref{Pf-10} that
\begin{align}\label{Pf-16}
\int_{B_R}&[f^{1-\frac{2}{p}}\nabla{f}+(p-2)
f^{{1-\frac{4}{p}}}(\nabla{f}\nabla{v})\nabla{v}]\nabla{\psi}+\beta_2(n,p,\sigma)\int_{B_R} f^2\psi\notag\\\leq&(n-1)Kp\int_{B_R} f^{2-\frac{2}{p}}\psi-\Big[\frac{2(p-1)}{n-1}-p\Big]\int_{B_R} f^{1-\frac{2}{p}}\psi\nabla{f}\nabla{v}.
\end{align}

In particular, if we denote
\begin{equation}\label{Pf-17}
\beta(n,p,\sigma)=\left\{\begin{array}{l}
\beta_1(n,p,\sigma), \quad a>0,\\
\beta_2(n,p,\sigma), \quad a<0,
\end{array}\right.
\end{equation}
then \eqref{Pf-13} together with \eqref{Pf-16} gives
\begin{align}\label{Pf-18}
\int_{B_R}&[f^{1-\frac{2}{p}}\nabla{f}+(p-2)
f^{{1-\frac{4}{p}}}(\nabla{f}\nabla{v})\nabla{v}]\nabla{\psi}+\beta(n,p,\sigma)\int_\Omega f^2\psi\notag\\
\leq&(n-1)Kp\int_\Omega f^{2-\frac{2}{p}}\psi-\Big[\frac{2(p-1)}{n-1}-p\Big]\int_{B_R} f^{1-\frac{2}{p}}\psi\nabla{f}\nabla{v}.
\end{align}

Let $\psi=f^b\eta^2$, where $\eta\in C_0^{\infty}(B_R)$ is non-negative and $ b>1$ is a constant to be determined. Then
\eqref{Pf-18} becomes
\begin{align}\label{Pf-19}
&\int_{B_R} bf^{b-\frac{2}{p}}|\nabla{f}|^2\eta^2+\int_{B_R} b(p-2)f^{{b-\frac{4}{p}}}
\eta^2(\nabla{f}\nabla{v})^2+2\int_{B_R} f^{b+1-{\frac{2}{p}}}\eta\nabla{f}\nabla{\eta} \notag\\
&+\int_{B_R} 2(p-2)f^{b+1-\frac{4}{p}}\eta(\nabla{f}\nabla{v})(\nabla{\eta}\nabla{v})+\beta(n,p,\sigma)\int_{B_R} f^{b+2}\eta^2\notag\\
\leq&(n-1)Kp\int_{B_R} f^{b+2-{\frac{2}{p}}}\eta^2-\Big[\frac{2(p-1)}{n-1}-p\Big]\int_{B_R} f^{b+1-{\frac{2}{p}}}\eta^2\nabla{f}\nabla{v},
\end{align}
which shows
\begin{align}\label{Pf-20}
&b\min\{1,p-1\}\int_{B_R} f^{b-\frac{2}{p}}|\nabla{f}|^2\eta^2+\beta(n,p,\sigma)\int_{B_R} f^{b+2}\eta^2\notag\\
\leq&(n-1)Kp\int_{B_R} f^{b+2-\frac{2}{p}}\eta^2-\Big[\frac{2(p-1)}{n-1}-p\Big]\int_{B_R} f^{b+1-\frac{2}{p}}\eta^2\nabla{f}\nabla{v}\notag\\
&-2\int_{B_R} f^{b+1-\frac{2}{p}}\eta\nabla{f}\nabla{\eta}-\int_{B_R}2(p-2)f^{b+1-\frac{4}{p}}\eta(\nabla{f}\nabla{v})(\nabla{\eta}\nabla{v})\notag\\
\leq&(n-1)Kp\int_{B_R} f^{b+2-{\frac{2}{p}}}\eta^2-\Big[\frac{2(p-1)}{n-1}-p\Big]\int_{B_R}  f^{b+1-{\frac{1}{p}}}|\nabla{f}|\eta^2\notag\\
&+2\big(\big|p-2\big|+1\big)\int_{B_R} f^{b+1-{\frac{2}{p}}}|\nabla{f}||\nabla{\eta}|\eta.
\end{align}
Using the Cauchy inequality
$$\aligned
2(|p-2|+1)f^{b+1-\frac{2}{p}}|\nabla{f}||\nabla{\eta}|\eta
\leq&\frac{b}{4}\min\{1,p-1\}f^{b-{\frac{2}{p}}}|\nabla{f}|^2\eta^2\notag\\
&+{\frac{4(|p-2|+1)^2}{b\min\{1,p-1\}}}f^{b+2-{\frac{2}{p}}}|\nabla{\eta}|^2
\endaligned$$
and
$$\aligned
&\Big[p-\frac{2(p-1)}{n-1}\Big]f^{b+1-{\frac{1}{p}}}|\nabla{f}|\eta^2\notag\\
&\leq\frac{b}{4}\min\{1,p-1\}f^{b-{\frac{2}{p}}}|\nabla{f}|^2\eta^2
+{\frac{\big[p-\frac{2(p-1)}{n-1}\big]^2}{b\min\{1,p-1\}}}f^{b+2}\eta^2,
\endaligned$$
then \eqref{Pf-20} becomes
\begin{align}\label{Pf-21}
&\int_{B_R}{\frac{b\min\{1,p-1\}}{2}}f^{b-{\frac{2}{p}}}|\nabla{f}|^2\eta^2
+\Big(\beta(n,p,\sigma)-{\frac{\big[p-\frac{2(p-1)}{n-1}\big]^2}{b\min\{1,p-1\}}}\Big)\int_{B_R} f^{b+2}\eta^2\notag\\
&\leq(n-1)Kp\int_{B_R} f^{b+2-{\frac{2}{p}}}\eta^2+{\frac{4\big(\big|p-2\big|+1\big)^2}{b\min\{1,p-1\}}}\int_{B_R} f^{b+2-{\frac{2}{p}}}|\nabla{\eta}|^2,
\end{align}
which gives
\begin{align}\label{Pf-22}
&\int_{B_R} \frac{b\min\{1,p-1\}}{2}f^{b-\frac{2}{p}}|\nabla{f}|^2\eta^2+\frac{\beta(n,p,\sigma)}{2}\int_{B_R} f^{b+2}\eta^2\notag\\
&\leq(n-1)Kp\int_{B_R} f^{b+2-\frac{2}{p}}\eta^2+\frac{4\big(\big|p-2\big|+1\big)^2}{b\min\{1,p-1\}}\int_{B_R} f^{b+2-\frac{2}{p}}|\nabla{\eta}|^2
\end{align}
by taking
$$b>\max\Big\{1,\frac{2\big[p-\frac{2(p-1)}{n-1}\big]^2}{\beta(n,p,\sigma)\min\{1,p-1\}}\Big\}$$
such that
$\frac{\beta(n,p,\sigma)}{2}>\frac{\big[p-\frac{2(p-1)}{n-1}\big]^2}{b\min\{1,p-1\}}$.
By the Cauchy inequlity
$$\aligned
|\nabla(f^{{\frac{b}{2}}+1-\frac{1}{p}}\eta)|^2\leq2\Big({\frac{b}{2}}+1
-\frac{1}{p}\Big)^2f^{b-\frac{2}{p}}|\nabla{f}|^2\eta^2+2f^{b+2-\frac{2}{p}}|\nabla{\eta}|^2
\endaligned$$
and multiplying both sides of the above inequality by a positive constant ${\frac{1}{2}}\big(\frac{b}{2}+1-{\frac{1}{p}}\big)^{-2}$, it follows from \eqref{Pf-22} that
\begin{align}\label{Pf-23}
&{\frac{b\min\{1,p-1\}}{4\big({\frac{b}{2}}+1-{\frac{1}{p}}\big)^2}}
\int_{B_R} |\nabla\big(f^{{\frac{b}{2}}+1-{\frac{1}{p}}}\eta\big)|^2
-\frac{b\min\{1,p-1\}}{2\big({\frac{b}{2}}+1-{\frac{1}{p}}\big)^2}\int_{B_R}
f^{b+2-{\frac{2}{p}}}|\nabla{\eta}|^2\notag\\
&+\frac{\beta(n,p,\sigma)}{2}\int_{B_R} f^{b+2}\eta^2\notag\\
\leq&(n-1)Kp
\int_{B_R} f^{b+2-{\frac{2}{p}}}\eta^2+{\frac{4\big(\big|p-2\big|+1\big)^2}{b\min\{1,p-1\}}}\int_{B_R} f^{b+2-{\frac{2}{p}}}|\nabla{\eta}|^2,
\end{align}
which is equivalent to
\begin{align}\label{Pf-24}
&\int_{B_R}|\nabla(f^{{\frac{b}{2}}+1-{\frac{1}{p}}}\eta)|^2
+\frac{\beta(n,p,\sigma)\big(b+2-{\frac{2}{p}}\big)^2}{2b\min\{1,p-1\}}\int_{B_R} f^{b+2}\eta^2\notag\\
\leq&\frac{Kp(n-1)\big(b+2-{\frac{2}{p}}\big)^2}{b\min\{1,p-1\}}\int_{B_R} f^{b+2-{\frac{2}{p}}}\eta^2\notag\\
&+\Big[{\frac{4\big(\big|p-2\big|+1\big)^2\big(b+2-\frac{2}{p}\big)^2}
{(b\min\{1,p-1\})^2}}+2\Big]\int_{B_R} f^{b+2-\frac{2}{p}}|\nabla\eta|^2.
\end{align}

The following Sobolev embedding theorem of Saloff-Coste plays a key role in our approach.

\begin{lem}\label{Le-2} (Theorem 3.1 in \cite{SC1992})
Let $M^n$ be a complete Riemannian manifold with $R_{ij}\geq-(n-1)K$. When $n>2$, for any $f\in C^\infty_0(B_R)$, there exists $C$ depending only on $n$, such that
\begin{align}\label{Pf-25}
\Big(\int_{B_R}& \big|f\big|^{2q}\Big)^\frac{1}{q}\leq\frac{e^{{C_0(n,p)}(1+\sqrt{K}R)}}{V^\frac{2}{n}}\Big[R^2\int_{B_R}|\nabla f|^2+\int_{B_R} f^2\Big],
\end{align}
where $q=\frac{n}{n-2}$ and $V$ denotes the volume of $B_R$. When $n=2$, the above inequality holds with $n$ replaced by any fixed $n'>2$.
\end{lem}

Therefore, the formula \ref{Pf-25} implies
\begin{align}\label{Pf-26}
\Big(\int_{B_R}& f^{\frac{2n}{n-2}\big({\frac{b}{2}}+1-\frac{1}{p}\big)}\eta^\frac{2n}{n-2}\Big)^\frac{n-2}{n}\notag\\
\leq&\frac{e^{{C_0(n,p)}(1+\sqrt{K}R)}}{V^{\frac{2}{n}}}\Big[R^2\int_{B_R}|\nabla(f^{{\frac{b}{2}}+1-\frac{1}{p}}\eta)|^2
+\int_{B_R}f^{b+2-\frac{2}{p}}\eta^2\Big].
\end{align}
According to \eqref{Pf-24} and \eqref{Pf-26}, we conclude that
\begin{align}\label{Pf-27}
&\Big(\int_{B_R} f^{\frac{2n}{n-2}\big(\frac{b}{2}+1-\frac{1}{p}\big)}\eta^\frac{2n}{n-2}\Big)^\frac{n-2}{2n}+{\frac{\beta(n,p,\sigma)\big(b+2-\frac{2}{p}\big)^2}{2nb\min\{1,p-1\}}}\frac{e^{{C_0(n,p)}(1+\sqrt{K}R)}}{V^\frac{2}{n}}R^2\int_{B_R}
f^{b+2}\eta^2\notag\\
\leq&\Big({\frac{Kp(n-1)\big(b+2-{\frac{2}{p}}\big)^2}{b\min\{1,p-1\}}}R^2+1\Big){\frac{e^{{C_0(n,p)}(1
+\sqrt{K}R)}}{V^\frac{2}{n}}}\int_{B_R} f^{b+2-\frac{2}{p}}\eta^2\notag\\&+\Big(\frac{4\big(\big|p-2\big|+1\big)^2\big(b+2
-{\frac{2}{p}}\big)^2}{(b\min\{1,p-1\})^2}+2\Big){\frac{e^{{C_0(n,p)}(1+\sqrt{K}R)}}{V^\frac{2}{n}}}R^2\int_{B_R} f^{b+2-\frac{2}{p}}|\nabla\eta|^2.
\end{align}
Let
$$b_0=C_1(n,p)(1+\sqrt{K}R),\ \ \ \ C_2=\frac{C_0(n,p)}{C_1(n,p)},$$
$$a_1={\frac{\beta(n,p,\sigma)\big(b+2-{\frac{2}{p}}\big)^2}{2nb^2\min\{1,p-1\}}},$$
$$a_2=\frac{1}{b^2_0b}\Big( {\frac{Kp(n-1)\big(b+2-{\frac{2}{p}}\big)^{2}}{b\min\{1,p-1\}}}R^2+1\Big),$$
$$a_3={\frac{4\big(\big|p-2\big|+1\big)^2\big(b+2-{\frac{2}{p}}\big)^2}{(b\min\{1,p-1\})^2}}+2.$$
Then, \eqref{Pf-27} becomes
\begin{align}\label{Pf-28}
&\Big(\int_{B_R} f^{\frac{2n}{n-2}\big(\frac{b}{2}+1-\frac{1}{p}\big)}\eta^\frac{2n}{n-2}\Big)^\frac{n-2}{2n}
+a_1b\frac{e^{c_2b_0}}{V^\frac{2}{n}}R^2\int_{B_R}  f^{b+2}\eta^2\notag\\
\leq&a_2b^2_0b\frac{e^{c_2b_0}}{V^\frac{2}{n}}\int_{B_R}  f^{b+2-\frac{2}{p}}\eta^2+a_3\frac{e^{c_2b_0}}{V^\frac{2}{n}}\int_{B_R}  f^{b+2-\frac{2}{p}}R^2|\nabla\eta|^2,
\end{align}
where $b_0>\max\{1,\frac{2\big[p-\frac{2(p-1)}{n-1}\big]^2}{\beta(n,p,\sigma)\min\{1,p-1\}}\}$.

\begin{lem}\label{Le-3}
Let $b_1=(b_0+2-\frac{2}{p}){\frac{n}{n-2}}$. Then there exist $C_3(n,p,\sigma)$, such that
\begin{align}\label{Pf-29}
\|f\|_{L^{b_1}(B_\frac{3R}{4})}\leq C_3(n,p,\sigma)\Big(\frac{b_0}{R}\Big)^pV^\frac{1}{b_1}.
\end{align}
\end{lem}

\proof We let
$$B_A=\Big\{x\in B_R;\ f>\big(\frac{2a_2}{a_1}\big)^\frac{p}{2}\big( \frac{b_0}{R}\big)^2\Big\},$$
$$B_B=\Big\{x\in B_R;\ f\leq\big(\frac{2a_2}{a_1}\big)^\frac{p}{2}\big(\frac{b_0}{R}\big)^2\Big\}.$$
If
$f>(\frac{2a_2}{a_1})^\frac{p}{2}(\frac{b_0}{R})^p$,
then
$$a_2b^3_0f^{b_0+2-\frac{2}{p}}<\frac{1}{2}a_1b_0R^2f^{b_0+2}.$$
If
$f\leq(\frac{2a_2}{a_1})^\frac{p}{2}(\frac{b_0}{R})^p$,
then
$$\int_{B_B}f^{b+2-\frac{2}{p}}\eta^2\leq\big(\frac{2a_2}{a_1}\big)^{\frac{p(b_0+2)}{2}-1}
\big(\frac{b_0}{R}\big)^{p(b_0+2)-2}\eta^2V.$$
Since
$B_R=B_A\cup B_B$,
we have
$$a_2b^2_0b\frac{e^{c_2b_0}}{V^\frac{2}{n}}\int_{B_A}f^{b+2-\frac{2}{p}}\eta^2
<\frac{1}{2}a_1b\frac{e^{c_2b_0}}{V^\frac{2}{n}}R^2\int_{B_R}f^{b+2}\eta^2,$$
$$a_2b^2_0b\frac{e^{c_2b_0}}{V^\frac{2}{n}}\int_{B_B}f^{b+2-\frac{2}{p}}\eta^2\leq
a_2b^3_0\frac{e^{c_2b_0}}{V^{\frac{2}{n}-1}}
\big(\frac{2a_2}{a_1}\big)^{\frac{p(b_0+2)}{2}-1}\big(\frac{b_0}{R}\big)^{p(b_0+2)-2}\eta^2.
$$
Therefore, letting $b=b_0$,
we obtain
\begin{align}\label{Pf-30}
&\Big(\int_{B_R} f^{\frac{2n}{n-1}(\frac{b_0}{2}+1-\frac{1}{p})}\eta^\frac{2n}{n-1}\Big)^\frac{n-1}{2n}
+\frac{1}{2}a_1b\frac{e^{c_2b_0}}{V^\frac{2}{n}}R^2\int_{B_R}f^{b_0+2}\eta^2\notag\\
&\leq a^{b_0}_4b^3_0\frac{e^{c_2b_0}}{V^{\frac{2}{n}-1}}\big(\frac{b_0}{R}\big)^{p(b_0+2)-2}
+a_3\frac{e^{c_2b_0}}{V^\frac{2}{n}}\int_{B_R}f^{b_0+2-\frac{2}{p}}R^2|\nabla\eta|^2,
\end{align}
where $a_4=\Big(a_2\big({\frac{2a_2}{a_1}}\big)^{\frac{p(b_0+2)}{2}-1}\Big)^{\frac{1}{b_0}}$.
Assume $0\leq\eta_1\leq1$, $\eta_1\equiv1$ in $B_{\frac{3R}{4}}$ and $|\nabla{\eta_1}|\leq{\frac{C_4(n)}{R}}$. Let $\eta=\eta^{b_0+2}_1$. By a direct calculation, we have
$$R^2|\nabla{\eta}|^2\leq(b_0+2)^2C_4(n)\eta^{\frac{2(b_0+1)}{b_0+2}}.$$
Denote by
$a_5=({\frac{b_0+2}{b_0}})^2C_4(n)$.
Thus
$$R^2|\nabla{\eta}|^2\leq a_5b^2_0\eta^{\frac{2(b_0+1)}{b_0+2}}.$$
Using the Yang inequality, the part of second term on the right hand of \eqref{Pf-30} can be written as
\begin{align}\label{Pf-31}
a_3\int_{B_R} f^{b_0+2-\frac{2}{p}}R^2|\nabla\eta|^2&\leq a_3a_5b^2_0\int_{B_R} f^{b_0+2-\frac{2}{p}}\eta ^\frac{2(b_0+1)}{b_0+2}\notag\\
&\leq a_3a_5b^2_0 \Big(\int_{B_R} f^{b_0+2}\eta^\frac{2(b_0+1)}{b_0+2-{\frac{2}{p}}}\Big)^{\frac{b_0+2
-\frac{2}{p}}{b_0+2}}V^{\frac{\frac{2}{p}}{b_0+2}}\notag\\
&\leq\frac{a_4}{2}b_0R^2\int_{B_R} f^{b_0+2}\eta^\frac{2(b_0+1)}{b_0+2-\frac{2}{p}}+a_6{\frac{b^{\frac{1}{2}p(b_0+2)+1}}{R^{p(b_0+2)-2}}}V\notag\\
&\leq\frac{a_4}{2}b_0R^2\int_{B_R} f^{b_0+2}\eta^\frac{2(b_0+1)}{b_0+2-\frac{2}{p}}+a_6b^3_0{\frac{b_0^{p(b_0+2)-2}}{R^{p(b_0+2)-2}}}V,
\end{align}
where
$$a_6=\frac{2(a_3a_5)^\frac{p(b_0+2)}{2}}{p(b_0+2)}\Bigg({\frac{2(b_0+2-{\frac{2}{p}}) }{a_4(b_0+2)}}\Bigg)^\frac{p(b_0+2)-2}{2}.$$
We arrive at
\begin{align}\label{Pf-32}
\Big(\int_{B_R}  f^{\frac{2n}{n-2}(\frac{b_0}{2}+1-\frac{1}{p})}\eta^\frac{2n}{n-2}\Big)^\frac{n-2}{2n}\leq (a^{b_0}_4+a_6)b^3_0\frac{e^{c_2b_0}}{V^{\frac{2}{n}-1}}\big(\frac{b_0}{R}\big)^{p(b_0+2)-2}.
\end{align}
Recall $b_1=(b_0+2-\frac{2}{p})\frac{n}{n-2}$. Taking the $\frac {1}{b_0+2-\frac{2}{p}}$ root on both side of \eqref{Pf-32}, we have
\begin{align}\label{Pf-33}
\|f\|_{L^{b_1}(B_{\frac{3R}{4}})}\leq a_7\big(\frac{b_0}{R}\big)^pV^\frac{1}{b_1},
\end{align}
where $a_7=\big[(a^{b_0}_4+a_6)b^3_0e^{c_2b_0}\big] ^\frac{1}{b_0+2-\frac{2}{p}}$. The proof of the lemma is completed.
\endproof

Let $a_8=\max\{a_2,a_3\}$. Then for \eqref{Pf-28}, by ignoring second term on the left hand, we have
\begin{align}\label{Pf-34}
&\Big(\int_{B_R} f^{{\frac{2n}{n-2}}({\frac{b}{2}}+1-{\frac{1}{p}})}\eta^{\frac{2n}{n-2}}\Big)^{{\frac{n-2}{2n}}}\leq a_8\frac{e^{c_2b_0}}{V^\frac{2}{n}}\int_{B_R}(b^2_0b\eta^2+R^2|\nabla{\eta}|^2)f^{b+2-\frac{2}{p}}.
\end{align}
To apply the Moser iteration, we set
$$b_{l+1}=b_l\frac{n}{n-2},\ B_l=B(O,\frac{R}{2}+\frac{R}{4^l}),\ l=1,2...$$
and choose $\eta_l$ such that $\eta_l\equiv1$ in $B_{l+1}$, $\eta_l\equiv0$ in $B_R\backslash B_l$,
$|\nabla \eta_l|\leq \frac{C_5{4^l}}{R},0\leq\eta_l\leq1$.
By letting $b+2-\frac{2}{p}=b_l,\eta=\eta_l$, we have
\begin{align}\label{Pf-35}
&\Big(\int_{B_{l+1}}f^{b_{l+1}}\Big)^\frac{1}{b_{l+1}}\leq  a_8^\frac{1}{b_{l}}\Big({\frac{e^{c_2b_0}}{V^\frac{2}{n}}}\Big)^\frac{1}{b_{l}}
\Big(\int_{B_l}(b^2_0b_l+R^2|\nabla{\eta_l}|^2)f^{b_l}\Big)^\frac{1}{b_l}.
\end{align}
Therefore,
\begin{align}\label{Pf-36}
\|f\|_{L^{{b_{l+1}}}({B_{l+1}})}\leq a_8^\frac{1}{b_{l}}\Big(\frac{e^{c_2b_0}}{V^\frac{2}{n}}\Big)^\frac{1}{b_{l}}(b^2_0b_l+C^2_5 16^l)^\frac{1}{b_l}\|f\|_{L^{{b_l}}({B_l})}.
\end{align}
Since $b_l=\big(b_0+2-\frac{2}{p}\big)\big(\frac{n}{n-2}\big)^{l}$ and $p>1$, we have
\begin{align}\label{Pf-37}
\|f\|_{L^{{b_{l+1}}}({B_{l+1}})}\leq& a_8^{\frac{1}{b_{l}}}\Big(\frac{e^{c_2b_0}}{V^\frac{2}{n}}\Big)^\frac{1}{b_{l}}
\Big({\frac{b_0+2-\frac{2}{p}}{b_0}}\Big)^\frac{1}{b_l}\notag\\
&\times\Big[b^3_0\Big(\frac{n}{n-2}\Big)^l+{\frac{b_0C^2_5 16^l}{b_0+2-\frac{2}{p}}}\Big]^\frac{1}{b_l}\|f\|_{L^{{b_l}}({B_l})}\notag\\
\leq& a_8^{\frac{1}{b_{l}}}\Big(\frac{e^{c_2b_0}}{V^\frac{2}{n}}\Big)^\frac{1}{b_{l}}
\Big({\frac{b_0+2-\frac{2}{p}}{b_0}}\Big)^\frac{1}{b_l}\notag\\
&\times\Big[b^3_0\Big(\frac{n}{n-2}\Big)^l+16^lC^2_5\Big]^\frac{1}{b_l}\|f\|_{L^{{b_l}}({B_l})}
\end{align}
Notice that $\sum^{\infty}_{l=1}\frac{1}{b_l}=\frac{n}{2b_1}$, $\sum^{\infty}_{l=1}\frac{l}{b_l}=\frac{n^2}{4b_1}$.
Using $(x+y)^\frac{1}{l}\leq x^\frac{1}{l}+y^\frac{1}{l}$ for any $x>0, y>0$ and positive integer $l$, we have
\begin{align}\label{Pf-38}
\prod^{\infty}_{l=1}\Big[b^3_0(\frac{n}{n-2})^l +16^lC^2_5\Big] ^\frac{1}{b_l}&=\prod^{\infty}_{l=1}\Big\{\Big[b^3_0\Big(\frac{n}{n-2}\Big)^l+16^lC^2_5\Big]^\frac{1}{l}\Big\}^\frac{l}{b_l}\notag\\
&\leq\prod^{\infty}_{l=1}(b^\frac{3}{l}_0\frac{n}{n-2}+16C^\frac{2}{l}_5)^\frac{l}{b_l}\leq\prod^{\infty}_{l=1}(C_6b^\frac{3}{l}_0)^\frac{l}{b_l}\notag\\&=C_6^{\sum^{\infty}_{l=1}\frac{l}{b_l}}(b^3_0)^{\sum^{\infty}_{l=1}\frac{1}{b_l}}=C_6^\frac{n^2}{4b_1}b_0^\frac{3n}{2b_1},
\end{align}
where constant $C_6$ is large enough such that
$C_6\geq\frac{n}{n-2}+16(\frac{C^2_5}{b^3_0})^\frac{1}{l}$.
Let us do a recurrence
\begin{align}\label{Pf-39}
\|f\|_{{L^{{\infty}}}(B_{\frac{R}{2}})}\leq&
a_8^{\sum^{\infty}_{l=1}{\frac{1}{b_l}}}\Big({\frac{e^{c_2b_0}}{V^\frac{2}{n}}}\Big)^{\sum^{\infty}_{l=1}{\frac{1}{b_l}}}
\Big({\frac{b_0+2-\frac{2}{p}}{b_0}}\Big)^{\sum^{\infty}_{l=1}{\frac{1}{b_l}}}\notag\\&\times\prod^{\infty}_{l=1}\Big[b^3_0\Big({\frac{n}{n-2}}\Big)^l+{\frac{b_0C^2_5 16^{l}}{b_0+2-{\frac{2}{p}}}}\Big]^{\frac{1}{b_l}}\|f\|_{{L^{{\infty}}}\big(B_{\frac{3R}{4}}\big)}\notag\\
\leq&a_8^{\frac{n}{2b_1}}\Big({\frac{e^{c_2b_0}}{V^\frac{2}{n}}}\Big)^{\frac{n}{2b_1}}
\Big({\frac{b_0+2-\frac{2}{p}}{b_0}}\Big)^{\frac{n}{2b_1}}C_6^\frac{n^2}{4b_1}b_0^\frac{3n}{2b_1}a_7\big(\frac{b_0}{R}\big)^pV^\frac{1}{b_1}\notag\\
=&a_7C_6^{{\frac{n^2}{4b_1}}}b_0^{{\frac{n}{b_1}}}\bigg(a_8\Big(b_0+2-{\frac{2}{p}}\Big)\bigg)^{\frac{n}{2b_1}}
e^{\frac{nc_2b_0}{2b_1}}\big({\frac{b_0}{R}}\big)^p.
\end{align}
Recalling that $b_0=C_1(n,p)(1+\sqrt{K}R)$ and $b_1=(b_0+2-\frac{2}{p})\frac{n}{n-2}$, using the fact that $\lim_{x\rightarrow \infty}x^\frac{1}{x}=1$ and $\lim_{R\rightarrow \infty}\frac{b_0}{b_1}=\frac{n-2}{n}$ over and over again, we obtain that
$$\sqrt{a_7C_1^pC_6^{{\frac{n^2}{4b_1}}}b_0^{{\frac{n}{b_1}}}
\bigg(a_8\Big(b_0+2-{\frac{2}{p}}\Big)\bigg)^{\frac{n}{2b_1}}e^{\frac{nc_2b_0}{2b_1}}}$$ convergence to a constant which depends only on $p,n$  and $\sigma$ if $R\rightarrow \infty$.
We conclude that
\begin{align}\label{Pf-40}
\frac{|\nabla u|}{u}\leq C_{n,p,\sigma}\Big(\frac{1+\sqrt{K}R}R\Big)^{\frac{p}{2}}
\end{align}
and then Theorem \ref{Th-1} follows.

\section{Proof of Theorem \ref{Th-2}}
Using the Cauchy inequality
$$\aligned
|\nabla v|^{2p-4}|{\rm Hess}\,v|^2_A\geq&\frac{1}{n}[|\nabla v|^{p-2}{\rm tr}_A({\rm Hess}\,v)]^2\notag\\
=&\frac{1}{n}(\Delta_{p}v)^2,
\endaligned$$
the Bochner formula \eqref{Pf-4} gives
\begin{align}\label{Res-2}
\frac{1}{p}\mathcal{L}(|\nabla v|^p)=&|\nabla v|^{2p-4}|{\rm Hess}\,v|^2_A+|\nabla v|^{2p-4}{\rm Ric}(\nabla v,\nabla v)\notag\\
&-|\nabla v|^{p-2}\nabla v\nabla[|\nabla v|^p+a(p-1)^{(p-1)}e^{(\frac{\sigma}{p-1}-1)v}]\notag\\\geq&\frac{1}{n}(\Delta_{p}v)^2-(n-1)K|\nabla v|^{2p-2}-|\nabla v|^{p-2}\nabla v\nabla |\nabla v|^{p}\notag\\
&-a(p-1)^{(p-1)}\Big(\frac{\sigma}{p-1}-1\Big)e^{(\frac{\sigma}{p-1}-1)v}|\nabla v|^{p}\notag\\
=&\frac{1}{n}|\nabla v|^{2p}+a(p-1)^{(p-1)}\Big(\frac{n+2}{n}-\frac{\sigma}{p-1}\Big)e^{(\frac{\sigma}{p-1}-1)v}|\nabla v|^{p}\notag\\
&-|\nabla v|^{p-2}\nabla v\nabla |\nabla v|^{p}-(n-1)K|\nabla v|^{2p-2}\notag\\
&+\frac{1}{n}a^2(p-1)^{2(p-1)}e^{2(\frac{\sigma}{p-1}-1)v}.
\end{align}
If two constants $a,\sigma$ satisfy that either $a>0$ and $\sigma\leq\frac{n+2}{n}(p-1)$, or $a<0$ and $\sigma\geq\frac{n+2}{n}(p-1)$, then
$$a\Big(\frac{n+2}{n}-\frac{\sigma}{p-1}\Big)\geq0.$$
Hence, from \eqref{Res-2}, we obtain
\begin{align}\label{Res-3}
\frac{1}{p}\mathcal{L}(|\nabla v|^p)\geq\frac{1}{n}|\nabla v|^{2p}-|\nabla v|^{p-2}\nabla v\nabla|\nabla v|^{p}-(n-1)K|\nabla v|^{2p-2},
\end{align}
which is equivalent to
\begin{align}\label{Res-4}
\mathcal{L}(f)\geq\frac{p}{n}f^2-(n-1)Kpf^{2-\frac{2}{p}}-pf^{1-\frac{2}{p}}\nabla{f}\nabla{v},
\end{align}
where $f=|\nabla v|^p$. Using the similar proof as in Theorem \ref{Th-1} by a suitable adjustment, we can complete the proof of Theorem \ref{Th-2}.

\bibliographystyle{Plain}

\begin{thebibliography}{10}

\bibitem{HZ2023}
G.Y. Huang, L. Zhao,
Liouville type theorems for $\Delta_pu+au^{\sigma}=0$ on complete noncompact Riemannian manifolds,
Chinese Ann. Math. Ser. B, {\em to appear}.

\bibitem{HHL2013}
G.Y. Huang, Z.J. Huang, H.Z. Li,
Gradient estimates and differential Harnack inequalities for a nonlinear parabolic equation on Riemannian manifolds,
Ann. Global Anal. Geom., 43(2013), 209-232.

\bibitem{2HHL2013}
G.Y. Huang, Z.J. Huang, H.Z. Li,
Gradient estimates for the porous medium equations on Riemannian manifolds,
J. Geom. Anal., 23(2013), 1851-1875.

\bibitem{HL2014}
G.Y. Huang, Z.J. Huang, H.Z. Li,
Gradient estimates and entropy formulae of porous medium and fast diffusion equations for the Witten Laplacian,
Pacific J. Math., 268(2014), 47-78.

\bibitem{HM2015}
G.Y. Huang, B.Q. Ma,
Gradient estimates and Liouville type theorems for a nonlinear elliptic equation,
Arch. Math.(Basel), 105(2015), 491-499.

\bibitem{PW2022}
P.L. Huang, Y.D. Wang,
Gradient estimates and Liouville theorems for Lichnerowicz equations,
Pacific J. Math., 317(2022), 363-386.

\bibitem{KN2009}
B. Kotschwar, L. Ni,
Local gradient estimates of $p$-harmonic functions, $1/H$-flow, and an entropy formula,
Ann. Sci. \'{E}c. Norm. Sup\'{e}r., 42(2009), 1-36.

\bibitem{Li1991}
J.Y. Li,
Gradient estimates and Harnack inequalities for nonlinear parabolic and nonlinear elliptic equations on Riemannian manifolds,
J. Funct. Anal., 100(1991), 233-256.

\bibitem{Ma2006}
L. Ma,
Gradient estimates for a simple elliptic equation on complete noncompact
Riemannian manifolds,
J. Funct. Anal., 241(2006), 374-382.

\bibitem{MH2017}
B.Q. Ma, G.Y.  Huang,
Hamilton-Souplet-Zhang's gradient estimates for two weighted nonlinear parabolic equations,
Appl. Math. J. Chinese Univ. Ser. B, 32(2017), 353-364.

\bibitem{MHL2018}
B.Q. Ma, G.Y.  Huang, Y. Luo,
Gradient estimates for a nonlinear elliptic equation on complete Riemannian manifolds,
Proc. Amer. Math. Soc., 146(2018),  4993-5002.

\bibitem{PWW2021}
B. Peng, Y.D. Wang, G.D. Wei,
Yau type gradient estimates for $\Delta u+au(\log u)^p+bu=0$ on Riemannian manifolds,
J. Math. Anal. Appl., 498(2021), Paper No. 124963, 24 pp.

\bibitem{SC1992}
L. Saloff-Coste, Laurent Uniformly elliptic operators on Riemannian manifolds,
J. Differential Geom., 36(1992), 417-450.

\bibitem{WZ2011}
X.D. Wang, L. Zhang,
Local gradient estimate for $p$-harmonic functions on Riemannian manifolds,
Comm. Anal. Geom., 19(2011), 759-771.

\bibitem{WW2023}
Y.D. Wang, G.D. Wei,
On the nonexistence of positive solution to $\Delta u+au^{p+1}=0$ on Riemannian manifolds,
J. Differential Equations, 362(2023), 7487.

\bibitem{Wang-Li-2016}
Y.Z. Wang, H.Q. Li,
Lower bound estimates for the first eigenvalue of the weighted $p$-Laplacian on smooth metric measure spaces,
Differential Geom. Appl., 45(2016), 23-42.

\bibitem{WY2021}
Y.Z. Wang, Y. Xue,
Gradient estimates for a parabolic $p$-Laplace equation with logarithmic nonlinearity on Riemannian manifolds,
Proc. Amer. Math. Soc., 149(2021), 1329-1341.

\bibitem{WY2023}
Y.Z. Wang, Y. Xue,
Logarithmic Harnack inequalities and differential Harnack estimates for $p$-Laplacian on Riemannian manifolds,
J. Math. Anal. Appl., 523(2023), Paper No. 127034.

\bibitem{Yang2008}
Y.Y. Yang,
Gradient estimates for a nonlinear parabolic equation on Riemannian manifolds,
Proc. Amer. Math. Soc., 136(2008), 4095-4102.

\bibitem{Yang2010}
Y.Y. Yang,
Gradient estimates for the equation $\Delta u+cu^{-\alpha}=0$ on Riemannian manifolds,
Acta. Math. Sinica, 26(2010), 1177-1182.

\bibitem{LY1975}
S.-T. Yau,
Harmonic functions on complete Riemannian manifolds,
Comm. Pure Appl. Math., 28(1975), 201-228.

\bibitem{ZY2018}
L. Zhao, D.Y. Yang,
Gradient estimates for the $p$-Laplacian Lichnerowicz equation on smooth metric measure spaces,
Proc. Amer. Math. Soc., 146(2018), 5451-5461.

\bibitem{ZS2020}
L. Zhao, M. Shen,
Gradient estimates for $p$-Laplacian Lichnerowicz equation on noncompact metric measure space,
Chinese Ann. Math. Ser., B 41(2020), 397-406.


\end{thebibliography}

\end{document}